\RequirePackage{fix-cm}
\documentclass[smallextended, natbib]{svjour3}       % onecolumn (second format)

\bibpunct[,]{[}{]}{,}{n}{,}{,}

\smartqed  % flush right qed marks, e.g. at end of proof

\usepackage{epsfig,amsmath,amssymb}
\usepackage{graphicx}
\usepackage{enumerate}
\usepackage{xfrac}
\usepackage{pgf}
\usepackage{pgfpages}
\usepackage{fix-cm}

\usepackage{caption}
\usepackage{subcaption}
\usepackage{dsfont}

\usepackage{amsfonts}
\usepackage{listings}

\usepackage[hidelinks]{hyperref}
\usepackage{booktabs} 
\usepackage{multirow}
\usepackage{epstopdf}

\usepackage{array}
\usepackage{setspace}
\usepackage{cancel}
\usepackage[noend]{algpseudocode}
\usepackage{algorithm}
\usepackage{algorithmicx}
\usepackage{mathrsfs}
\usepackage{upgreek}
\usepackage{eufrak}
\usepackage{enumerate}
\usepackage{accents}
\usepackage{arydshln}
\usepackage{textgreek}
\usepackage{multicol}
\usepackage{lscape}
\usepackage[toc,page]{appendix}
\usepackage{blindtext, rotating}
\usepackage{moresize}
\usepackage{scalerel}
\usepackage[english]{babel}
\usepackage{tabularx}

\DeclareMathAlphabet{\mathpzc}{OT1}{pzc}{m}{it}

\DeclareMathAlphabet{\mathcal}{OMS}{cmsy}{m}{n}

\makeatletter
\newcommand{\pushright}[1]{\ifmeasuring@#1\else\omit\hfill$\displaystyle#1$\fi\ignorespaces}
\newcommand{\pushleft}[1]{\ifmeasuring@#1\else\omit$\displaystyle#1$\hfill\fi\ignorespaces}
\makeatother

\def\R{\hbox{$\mathbb R$}}

\renewcommand{\bf}{\mathbf}

\newcommand{\bul}{\bullet}
\newcommand{\affineset}{\mathscr{A}}
\newcommand{\paramspace}{\Theta}

\newcommand{\B}{\mathbin{\ThisStyle{\scalebox{1.15}{$\SavedStyle\ensuremath{\mathpzc{B}}$}}}}

\spnewtheorem{assumption}{Assumption}{\bfseries}{\itshape}
\spnewtheorem{observation}{Observation}{\bfseries}{\itshape}

\def\dataname{Data Availability}
\def\data{\par\addvspace{17pt}\small\rmfamily
\trivlist\if!\dataname!\item[]\else
\item[\hskip\labelsep
{\bfseries\dataname}]\fi}

\def\decname{Declarations}
\def\declarations{\par\addvspace{17pt}\small\rmfamily
\trivlist\if!\decname!\item[]\else
\item[\hskip\labelsep
{\bfseries\decname}]\fi}

\begin{document}

\title{A Unified Tool for Solving Uni-Parametric Linear Programs, Convex Quadratic Programs, and Linear Complementarity Problems}
\titlerunning{A Unified Tool for Solving upLP, upQP, and upLCP}

\author{Nathan Adelgren}

\institute{Nathan Adelgren \at
              Andlinger Center for Energy and the Environment\\
              Princeton University\\
              Princeton, NJ, 08544\\
              \email{na4592@princeton.edu}, ORCHID: \url{https://orcid.org/
0000-0003-3836-9324}           %  \\
}

\date{Received: date / Accepted: date}

\maketitle

\begin{abstract}
We introduce a new technique for solving uni-parametric versions of linear programs, convex quadratic programs, and linear complementarity problems in which a single parameter is permitted to be present in any of the input data. We demonstrate the use of our method on a small, motivating example and present the results of a small number of computational tests demonstrating its utility for larger scale problems.
\keywords{parametric optimization \and linear complementarity \and linear programming \and quadratic programming}
%\subclass{MSC code1 \and MSC code2 \and more}
\end{abstract}

\begin{acknowledgements}
The author would like to thank Jacob Adelgren for offering helpful advice and feedback regarding the implementation.
\end{acknowledgements}

\section{Introduction}\label{sec:intro}

In this work we consider the uni-parametric form of the Linear Complementarity Problem (LCP) in which all input data is permitted to be dependent on a single parameter $\theta \in \paramspace$, where 
\begin{equation}\label{eq:paramspace}
\paramspace := \{\theta \in \R: \alpha \leq \theta \leq \beta \text{ and } \alpha,\beta \in \R\}
\end{equation} 
is a connected interval in $\R$ that represents the set of ``attainable'' values for $\theta$. This problem is referred to as the uni-parametric (or single-parametric) Linear Complementarity Problem (upLCP). Let $\affineset = \{\mu \theta + \sigma: \mu, \sigma \in \R\}$, the set of affine functions of $\theta$. Then upLCP is as follows:
\begin{quote}
Given $M(\theta)\in \affineset^{h\times h}$ and $q(\theta)\in \affineset^h$, for each $\theta\in \paramspace$ find vectors $w(\theta)$ and $z(\theta)$ that satisfy the system
\begin{equation}\label{upLCP}
\begin{array}{c}
w - M(\theta)z = q(\theta)\\[1mm]
w^\top z = 0\\[1mm]
w,z \geq 0
\end{array}
\end{equation}
or show that no such vectors exist.
\end{quote}

upLCP is said to be \emph{feasible at $\theta$} if there exist $w(\theta)$ and $z(\theta)$ that satisfy System \eqref{upLCP} and \emph{infeasible at $\theta$} otherwise. Similarly, upLCP is said to be \emph{feasible} if there exists a $\hat{\theta} \in \paramspace$ at which upLCP is feasible, and \emph{infeasible} otherwise. 

Now, recognize that $\paramspace$ must be an infinite set, otherwise upLCP reduces to LCP. Hence, it is not possible to determine a solution to System \eqref{upLCP} for each $\theta \in  \paramspace$ individually. Instead, upLCP is solved by partitioning the interval $\paramspace$ into a set of \emph{invariancy intervals}. As the name ``invariancy intervals'' suggests, within each of these intervals the representation of the solution vectors $w$ and $z$ as functions of $\theta$ is invariant.  The methods we propose solve upLCP whenever the following assumptions are met.

\begin{assumption}\label{asm:sufficient}
The matrix $M(\theta)$ is sufficient for all $\theta \in \paramspace$.
\end{assumption}

\begin{assumption}\label{asm:feasible}
System \eqref{upLCP} is feasible for all $\theta \in \paramspace$.
\end{assumption}

We point any reader interested in background information on LCP, in particular the definition of sufficient matrices, to the work of \citet{cottle2009linear}. Additionally, we note that Assumption \ref{asm:feasible} is not extremely restrictive. For example, the claim of Assumption \ref{asm:feasible} is satisfied when upLCP results from the reformulation of a biobjective linear program or biobjective quadratic program that has been scalarized using the weighted sum approach (see \citep{ehrgott2005multicriteria}, for example) and $\paramspace$ is taken to be $[0,1]$.

To our knowledge, the only other works in which solution procedures are proposed for upLCP having the general form of System \eqref{upLCP} are those of \citet{valiaho1994procedure} and \citet{chakraborty2004solution}. Section 6 of the former work describes various restrictions that can be placed on the structure of $M(\theta)$ in order to guarantee finite execution of the procedures presented therein, but all are more restrictive than our Assumption \ref{asm:sufficient}. Similarly, the methods presented in the latter work require that $M(\theta)$ be a $P$-matrix for all $\theta \in \paramspace$, which is again more restrictive than our Assumption \ref{asm:sufficient}. Solution methodology for the multiparametric counterpart to System \eqref{upLCP}, i.e., the case in which $\theta \in \R^k$ for $k > 1$, is presented in \citep{adelgren2021advancing}. However, the techniques presented therein are not directly applicable for scalar $\theta$. This motivates our current work, the remainder of which is organized as follows. The proposed algorithm for solving upLCP is presented in Section \ref{sec:algorithm}. We demonstrate the use of our proposed methodology on a motivating example in Section \ref{sec:example}. In Section \ref{sec:otherProblems} we discuss some classes of optimization problems that can be reformulated as upLCP and can therefore be solved using the methods presented in this work. %Section \ref{sec:implementation/results} contains a brief description of a software implementation of our proposed methods that we make available to the research community as well as the results of a small number of computational tests.
Section \ref{sec:implementation/results} contains a brief description of our implementation and the results of a small number of computational tests.
%are presented in Section \ref{sec:results}. 
Finally, we give concluding remarks in Section \ref{sec:conclusion}.

\section{Proposed Algorithm}\label{sec:algorithm}

We begin this section with a brief set of definitions. Given an instance of upLCP, as presented in System \eqref{upLCP}, we define the matrix $G(\theta):= \left[ \begin{matrix}
I & -M(\theta)
\end{matrix} \right]$ and the index set $\mathcal{E} := \{1,\dots,2h\}$. Then a set $\B \subset \mathcal{E}$ is called a \emph{basis} if $|\B|=h$. Moreover, a basis $\B$ is \emph{complementary} if $\left|\{i,i+h\}\cap \B \right| = 1$ for each $i \in \{1,\dots,h\}$. Then, given a complementary basis $\B$, its associated \emph{invariancy interval} is the set 
\begin{equation}\label{invInterval} \mathcal{II}_{\B} := \left\{ \theta \in \paramspace: G(\theta)^{-1}_{\bul \B} q(\theta) \geq 0 \right\},
\end{equation}
where $G(\theta)^{-1}_{\bul \B}$ denotes the matrix comprised of the columns of $G(\theta)^{-1}$ whose indices are in $\B$. We note here that for a given complementary basis $\B$, $\mathcal{II}_{\B}$ may be the union of disjoint intervals in $\paramspace$. To see this, recognize that: (i) for any complementary basis $\B$, the system $G(\theta)^{-1}_{\bul \B} q(\theta) \geq 0$ can be equivalently written as $\dfrac{Adj\left(G(\theta)_{\bul \B}\right)}{det\left(G(\theta)_{\bul \B}\right)} q(\theta) \geq 0$, where $Adj(\cdot)$ and $det(\cdot)$ represent the matrix adjoint and determinant, respectively; and (ii) both the adjoint and determinant of a given matrix can be represented as polynomials of the elements of the matrix. Hence, for each complementary basis $\B$, $\mathcal{II}_{\B}$ is a possibly nonconvex subset of $\paramspace$ defined by a set of rational inequalities in $\theta$. Fortunately, this structure can be somewhat improved. By Lemma 2.1 of \citep{adelgren2021advancing} we know that, under Assumption \ref{asm:sufficient}, the sign of $det\left(G(\theta)_{\bul \B}\right)$ is invariant over $\paramspace$ for any complementary basis $\B$. Thus, $\mathcal{II}_{\B}$ can be represented as 
\begin{equation}\label{invInterval2} \mathcal{II}_{\B} := \left\{ \theta \in \paramspace: s_{\B} Adj(G(\theta)_{\bul \B}) q(\theta) \geq 0 \right\},
\end{equation}
where $s_{\B}$ represents the sign of $det\left(G(\theta)_{\bul \B}\right)$ over $\paramspace$. Note that in Equation \ref{invInterval2}, $\mathcal{II}_{\B}$ is given by a system of inequalities that are polynomial in $\theta$. We now provide Algorithm \ref{alg:solveUpLCP} in which we present a technique for solving upLCP.
\begin{algorithm}%[H]
\small
  \caption{\small\textsc{Solve\_upLCP}($\paramspace$)~--~Partition the parameter space $\paramspace$.\\
 \textbf{Input}: The set $\paramspace$ as defined in Equation \eqref{eq:paramspace} for an instance of upLCP.\\
 \textbf{Output}: A set $\mathcal{P}$ specifying a partition of $\paramspace$. Each $\mathcal{I} \in \mathcal{P}$ is a tuple of the form $(\protect\B^*, \alpha^*, \beta^*)$ where $\protect\B^*$ represents a complementary basis and $\alpha^*,\beta^* \in \protect\R$ specify an interval $[\alpha^*,\beta^*] \subseteq \mathcal{II}_{\protect\B}$ such that $w(\theta)$ and $z(\theta)$ satisfy System \eqref{upLCP} for all $\theta \in [\alpha^*,\beta^*]$.
  }\label{alg:solveUpLCP}
\begin{algorithmic}[1]
\State Let $\mathcal{S} = \{\paramspace\}$ and $\mathcal{P} = \emptyset$.
\While{$\mathcal{S} \neq \emptyset$}{ select $[\alpha', \beta']$ from $\mathcal{S}$.}\label{alg:solveUpLCP:line:while}
	\State Set $\theta^* = \frac{\alpha' + \beta'}{2}$ and compute a complemetary basis $\B^*$ such that $\theta^* \in \mathcal{II}_{\B^*}$.\label{alg:solveUpLCP:line:midpoint}
	%\State Compute $\alpha^*, \beta^* \in \R$ such that $\theta^* \in [\alpha^*,\beta^*] \subseteq \mathcal{II}_{\B^*}$ and add $(\B^*, \alpha^*,\beta^*)$ to $\mathcal{P}$.\label{alg:solveUpLCP:line:endpoints}
	\State Set $\alpha^*,\beta^* = \textsc{Get\_Extremes}(\B^*, \theta^*, \alpha', \beta')$ and add $(\B^*, \alpha^*,\beta^*)$ to $\mathcal{P}$.\label{alg:solveUpLCP:line:endpoints}
%	\State Add $(\B^*, \alpha^*,\beta^*)$ to $\mathcal{P}$.
	\If{$\alpha' < \alpha^*$}{ add $[\alpha',\alpha^*]$ to $\mathcal{S}$.}\label{alg:solveUpLCP:line:subint1}
	\EndIf
	\If{$\beta' > \beta^*$}{ add $[\beta^*,\beta']$ to $\mathcal{S}$.}\label{alg:solveUpLCP:line:endwhile}
	\EndIf
\EndWhile
\State Return $\mathcal{P}$.
\end{algorithmic}
\end{algorithm}
The majority of the work done in Algorithm \ref{alg:solveUpLCP} is the processing of a set $\mathcal{S}$ of intervals using the while loop contained in lines \ref{alg:solveUpLCP:line:while}--\ref{alg:solveUpLCP:line:endwhile}. On line \ref{alg:solveUpLCP:line:midpoint}, the midpoint $\theta^*$ of the current interval is computed and a complementary basis $\B^*$ is sought for which $\theta^*$ is contained within $\mathcal{II}_{\B^*}$. Such a complementary basis can be found by fixing $\theta$ to $\theta^*$ and solving the resulting non-parametric LCP using, for example, the criss-cross method of \citet{den1993linear}. We note that the criss-cross method is a particularly good choice in this case as it is guaranteed to solve the resulting non-parametric LCP under Assumption \ref{asm:sufficient}. On line \ref{alg:solveUpLCP:line:endpoints}, the subroutine \textsc{Get\_Extremes} is used to compute $\alpha^*,\beta^* \in \R$ that define a subinterval $[\alpha^*,\beta^*]$ of $\mathcal{II}_{\B^*}$ that contains $\theta^*$ and should be included in the final partition of $\paramspace$. 
Finally, lines \ref{alg:solveUpLCP:line:subint1}--\ref{alg:solveUpLCP:line:endwhile} identify connected subintervals of $[\alpha',\beta']$ over which solutions to System \eqref{upLCP} have yet to be computed, if any exist, and those subintervals are added to the set $\mathcal{S}$.

We pause now to discuss the subroutine \textsc{Get\_Extremes} in more detail. There are many ways to implement such a routine, but the strategies employed in our implementation are presented in Algorithm \ref{alg:getExtremes}.
\begin{algorithm}%[H]
\small
  \caption{\small\textsc{Get\_Extremes}($\protect\B^*,\theta^*,\alpha',\beta'$)~--~Compute extreme values of $\mathcal{II}_{\protect\B}$.\\
 \textbf{Input}: A complementary basis $\protect\B^*$, $\theta^* \in \mathcal{II}_{\protect\B^*}$, and $\alpha',\beta' \in \R$.\\
 \textbf{Output}: $\alpha^*,\beta^* \in \protect\R$ that specify an interval $[\alpha^*,\beta^*] \subseteq \mathcal{II}_{\protect\B^*}$ such that $\theta^* \in [\alpha^*,\beta^*]$ and $w(\theta)$ and $z(\theta)$ satisfy System \eqref{upLCP} for all $\theta \in [\alpha^*,\beta^*]$.
  }\label{alg:getExtremes}
\begin{algorithmic}[1]
\State Set $\alpha^* = \alpha'$ and $\beta^* = \beta'$.
\For{$i \in \B^*$}\label{alg:getExtremes:lineFor}
	\If{$degree\left(\left(Adj(G(\theta)_{\bul \B})\right)_{i\bul} q(\theta)\right) > 0$}
		\State Let $\mathcal{R} = \textsc{Get\_Real\_Roots}\left(\left(Adj(G(\theta)_{\bul \B})\right)_{i\bul} q(\theta)\right)$.\label{alg:getExtremes:lineGetRoots}
		\For{$r \in \mathcal{R}$}\label{alg:getExtremes:lineFor2}
			\If{$multiplicity(r)$ is odd}{}\label{alg:getExtremes:lineMult}
				\If{$r \in (\alpha^*,\theta^*)$}{ set $\alpha^* = r$.}\label{alg:getExtremes:lineLeft}
				\ElsIf{$r \in (\theta^*,\beta^*)$}{ set $\beta^* = r$.}\label{alg:getExtremes:lineRight}
				\ElsIf{$r == \alpha^*$ OR $r == \beta^*$}{}\label{alg:getExtremes:lineEqual}
					\If{$\left.\frac{d}{d\theta}\left(s_{\B} \left(Adj(G(\theta)_{\bul \B})\right)_{i\bul} q(\theta) \right)\right|_{\theta = \theta^*} > 0$}{ set $\alpha^* = r$.}
					\Else{ set $\beta^* = r$.}\label{alg:getExtremes:lineRight2}
					\EndIf
				\EndIf
			\EndIf
		\EndFor
	\EndIf
\EndFor
\State Return $\alpha^*,\beta^*$.
\end{algorithmic}
\end{algorithm}
On a high-level, the work done in Algorithm \ref{alg:getExtremes} seeks to find the largest connected interval that is a subset of $\mathcal{II}_{\B^*}$ and contains $\theta^*$. To do this, we examine the roots of the nonconstant polynomial functions that serve as the boundaries of $\mathcal{II}_{\B^*}$ (lines \ref{alg:getExtremes:lineFor}--\ref{alg:getExtremes:lineGetRoots}). We note that any polynomial root finder capable of determining root multiplicity can be used for the subroutine \textsc{Get\_Real\_Roots} (see, for example, \citep{rouillier2004efficient}). As roots of even multiplicity cannot be restrictive, we need only consider roots of odd multiplicity (lines \ref{alg:getExtremes:lineFor2}--\ref{alg:getExtremes:lineMult}). Since we have $\theta^* \in \mathcal{II}_{\B^*}$, we know that if no root occurs directly at $\theta^*$, then the endpoints of the interval of interest are: (i) the root closest to $\theta^*$ on the left, and (ii) the root closest to $\theta^*$ on the right (lines  \ref{alg:getExtremes:lineLeft}--\ref{alg:getExtremes:lineRight}). If a root does occur directly at $\theta^*$, however, then $\theta^*$ serves as one of the endpoints of the interval of interest, and we must determine the direction in which the function associated with that root increases in order to determine which endpoint it is (lines \ref{alg:getExtremes:lineEqual}--\ref{alg:getExtremes:lineRight2}). For this purpose, we use the sign of the derivative of the associated function, evaluated at $\theta^*$. Since the defining constraints of $\mathcal{II}_{\B^*}$ are given as greater-than-or-equal-to constraints, we know that if the aforementioned sign is positive, $\theta^*$ is the left endpoint of the interval of interest. Otherwise, it is the right endpoint. We note that it is possible that two functions defining boundaries of $\mathcal{II}_{\B^*}$ each have roots at $\theta^*$ and, moreover, the derivatives of the two functions, evaluated at $\theta^*$, may have opposite signs. In this case, the interval of interest reduces to a singleton at $\theta^*$ and need not be included in the partition of $\paramspace$. Hence, such intervals can be rejected upon discovery or removed from the partition in a post-processing phase.

We finish this section by stating a small number of theoretical results related to the methods proposed in Algorithms \ref{alg:solveUpLCP}--\ref{alg:getExtremes}. For the sake of space, we do not include proofs for these results -- though, the proofs are not challenging -- and, as such, we state them as observations.

\begin{observation}\label{obs:multiple}
For a given complementary basis $\B$, the final partition $\mathcal{P}$ of $\paramspace$ may contain more than one connected subinterval of $\mathcal{II}_{\B}$.
\end{observation}

\begin{observation}
Under Assumption \ref{asm:sufficient}, for a given complementary basis $\B$, the number of connected subintervals of $\mathcal{II}_{\B}$ present in the final partition $\mathcal{P}$ of $\paramspace$ can be at most $n-h$, where $n$ is the number of unique roots of odd multiplicity of the functions $Adj(G(\theta)_{\bul \B}) q(\theta)$ and $h = |\B|$.
\end{observation}

\begin{observation}
Under Assumptions \ref{asm:sufficient}--\ref{asm:feasible}, Algorithm \ref{alg:solveUpLCP} will complete in a finite number of iterations.
\end{observation}

\section{Motivating Example}\label{sec:example}

\begin{figure}[t]
\caption{Regions over which the computed parametric solutions are valid.}\label{fig:ex2}
\centering
\includegraphics[width=70mm]{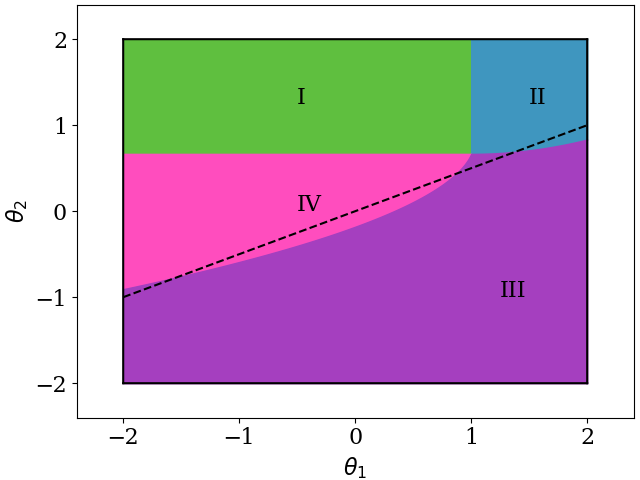} %.6\textwidth]{Figure1}
\end{figure}

At the start of Chapter 1 of \citep{murty1997linear}, the authors demonstrate that an instance of LCP with $M= \left[\begin{array}{rr}
2 & -1\\
1 & 3\\
\end{array} \right]$ will have relatively nice properties. We use this as a starting point and consider the following instance of mpLCP:
\begin{equation}\label{eq:ex2}
\begin{array}{c}
w - \left[\begin{array}{cc}
2 & -1 + \theta_2\\
1 - \theta_1 & 3\\
\end{array} \right]z = \left[\begin{array}{c}
1 - \theta_1\\
-2 + 3\theta_2\\
\end{array} \right]\\[1mm]
w^\top z = 0\\[1mm]
w,z \geq 0
\end{array}
\end{equation}
We assume here that $(\theta_1, \theta_2) \in [-2,2]^2$ and note that it is straightforward to verify that $M(\theta) = \left[\begin{array}{cc}
2 & -1 + \theta_2\\
1 - \theta_1 & 3\\
\end{array} \right]$ is sufficient for all $(\theta_1, \theta_2) \in [-2,2]^2$ using Theorem 4.3 of \citep{valiaho1996criteria}. We note that beginning with a mpLCP rather than a upLCP allows us to visualize some of the challenges associated with upLCP. We compute the parametric solution to System \eqref{eq:ex2} using the techniques proposed in \citep{adelgren2021advancing}. The final partition of $\paramspace$ is depicted in Figure \ref{fig:ex2} and the complementary bases associated with each region, as well as the parametric solutions associated with each basis, are provided in Table \ref{tab:ex2}. 
\begin{table} 
\caption{Parametric solutions to System \eqref{eq:ex2}.}\label{tab:ex2}
\centering
\begin{tabular}{ll}
\hline
Region 	& Basic Variables and\\
		& Associated Values\\
\hline
\hline
I  	& $w_1 = 1 - \theta_1$\\
	& $w_2 = 3\theta_2 - 2$\\
\hline
II	& $w_1 = -\theta_1 - \theta_2^2 + \frac{5}{3}\theta_2 + \frac{1}{3}$\\
	& $z_2 = \frac{2}{3} - \theta_2$\\
\hline
III	& $z_1 = \frac{3\theta_1 + 3\theta_2^2 - 5\theta_2 - 1}{\theta_1\theta_2 - \theta_1 - \theta_2 + 7}$\\
	& $z_2 = \frac{\theta_1^2 - 2\theta_1 - 6\theta_2 + 5}{\theta_1\theta_2 - \theta_1 - \theta_2 + 7}$\\
\hline
IV	& $z_1 = \frac{1}{2}\theta_1 - \frac{1}{2}$\\
	& $w_2 = -\frac{1}{2}\theta_1^2 + \theta_1 + 3\theta_2 - \frac{5}{2}$\\
\hline
\end{tabular}
\end{table}
We note that, in an abuse of notation, we identify a basis in Table \ref{tab:ex2} by its associated basic variables, rather than the elements of $\mathcal{E}$ that it contains.  Additionally, variables having zero value are omitted from Table \ref{tab:ex2}. The same convention for identifying complementary bases and omitting variables with zero value is used throughout the remainder of this work.

We now reduce System \eqref{eq:ex2} to an instance of upLCP by fixing $\theta_2 = \frac{1}{2}\theta_1$ (depicted in Figure \ref{fig:ex2} as a dotted black line) and letting $\paramspace = [-2,2]$. In observing Figure \ref{fig:ex2}, we note that this reduction will result in a realization of the scenario outlined in Observation \ref{obs:multiple}. Namely, our final partition of $\paramspace$ will contain more than one interval associated with a single complementary basis. We believe that situations like this have been one of the primary barriers to the development of methods for solving upLCP as specified by System \eqref{upLCP}. We now demonstrate how the techniques outlined in Algorithms \ref{alg:solveUpLCP}--\ref{alg:getExtremes} can be used to solve this instance of upLCP. We begin by setting $\mathcal{S} = [-2,2]$ and then iterate through the while loop on line \ref{alg:solveUpLCP:line:while} of Algorithm \ref{alg:solveUpLCP}. Note that in the following discussion we process elements from the set $\mathcal{S}$ using a last-in-first-out strategy.\\

\noindent \textbf{Iteration 1}:
\begin{quote}
Here $[\alpha',\beta'] = [-2,2]$ and we set $\theta^* = 0$. Using the criss-cross method, we find that an optimal basis at $\theta^*$ is $\B_1 =  \{w_1 = -\frac{1}{4}\theta_1^2 - \frac{1}{6}\theta_1 + \frac{1}{3}, z_2 = -\frac{1}{2}\theta_1 + \frac{2}{3}\}$. We now enter Algorithm \ref{alg:getExtremes} via the routine \textsc{Get\_Extremes}. The real roots associated with the function representing $w_1$ are approximately $-1.535$ and $0.869$, whereas the function associated with $z_2$ has only one real root, approximately $1.333$. Of these roots, $-1.535$ is closest to $\theta^*$ on the left and $0.869$ is closest on the right. As such, these values are returned from \textsc{Get\_Extremes} as $\alpha^*$ and $\beta^*$, respectively. On lines \ref{alg:solveUpLCP:line:subint1}--\ref{alg:solveUpLCP:line:endwhile} we now add $[-2,-1.535]$ and $[0.869,2]$ to $\mathcal{S}$ and proceed to Iteration 2.
\end{quote}

\noindent \textbf{Iteration 2}:
\begin{quote}
For this and subsequent iterations, we omit details that are analogous to those of Iteration 1. We note that here we have $[\alpha',\beta'] = [0.869,2]$ and thus $\theta^* = 1.4345$. An optimal basis at this point is $\B_2 = \{z_1 = \frac{1}{2}\theta_1 - \frac{1}{2}, w_2 = -\frac{1}{2}\theta_1^2 + \frac{5}{2}\theta_1 - \frac{5}{2} \}$ and respective sets of approximate real roots for these functions are $\{1\}$ and $\{1.382,3.618\}$. Hence, we set $\alpha^* = 1.382$ and leave $\beta^*$ as $2$. The only new interval added to $\mathcal{S}$ is then $[0.869,1.382]$.
\end{quote}

\noindent \textbf{Iteration 3}:
\begin{quote}
We now have $[\alpha',\beta'] = [0.869,1.382]$ and $\theta^* = 1.1255$. An optimal basis here is $\B_3 = \{z_1 = \frac{-3\theta_1^2 - 2\theta_1 + 4}{-2\theta_1^2 + 6\theta_1 - 28}, z_2 = \frac{-2\theta_1^2 + 10\theta_1 - 10}{-\theta_1^2 + 3\theta_1 - 14} \}$ and respective sets of approximate real roots for these functions are $\{0.869,-1.535\}$ and $\{1.382,3.618\}$. As a result, $\alpha^*$ and $\beta^*$ are left as $0.869$ and $1.382$, respectively, and no intervals are added to $\mathcal{S}$.
\end{quote}

\noindent \textbf{Iteration 4}:
\begin{quote}
Finally, we have $[\alpha',\beta'] = [-2,-1.535]$ and $\theta^* = -1.7675$. An optimal basis at this point is $\B_3$. Thus, using the previously computed sets of approximate real roots, we see that $\alpha^*$ and $\beta^*$ are left as $-2$ and $-1.535$, respectively, and no intervals are added to $\mathcal{S}$.
\end{quote}

\noindent We have now completed the execution of Algorithm \ref{alg:solveUpLCP}. The final partition of $\paramspace$ consists of the four intervals $[-1.535,0.869]$, $[1.382,2]$, $[0.869,1.382]$, and $[-2,-1.535]$. The parametric solutions that are valid over these intervals are those given in the respective descriptions of $\B_1, \B_2, \B_3$ and (again) $\B_3$ above. It is straightforward to verify that this solution matches the one given in Table \ref{tab:ex2} when $\theta_2$ is fixed to $\frac{1}{2}\theta_1$. Recognize, though, that the representation of $w_2$ in the solution associated with Region I in Table \ref{tab:ex2} reduces to $w_2 = \frac{3}{2}\theta_1 - 2$. Hence, as there exists no $\theta_1 \in [-2,2]$ for which both $w_1 = 1-\theta_1$ and $w_2 = \frac{3}{2}\theta_1 - 2$ are non-negative, the complementary basis $\B' = \{w_1,w_2\}$ is not associated with any interval in the final partition of $\paramspace$ when solving the reduced instance of upLCP. 

\section{Applicability to Other Classes of Problems}\label{sec:otherProblems}

It is well known that LCP arises naturally as the system resulting from applying the Karush-Kuhn-Tucker (KKT) optimality conditions to quadratic programs (QP). Moreover, as a solution to the KKT system is guaranteed to give an optimal solution for any convex QP, the methods presented herein are directly applicable to the uni-parametric form of QP given by
\begin{align}
\min_x\quad & \frac{1}{2} x^\top Q(\theta) x + c(\theta)^\top x\nonumber\\
\text{s.t.}\quad & A(\theta) x \leq b(\theta)\label{prob:QP}\\
& x \geq 0\nonumber\\
& \theta \in \paramspace\nonumber
\end{align}
whenever $Q(\theta)$ is convex for all $\theta \in \paramspace$, with $\paramspace$ defined as in Equation \eqref{eq:paramspace}. We note that this implies that the methods presented herein are also applicable to: (i) biobjective problems having linear and/or convex quadratic objectives and linear constraints (via weighted-sum scalarization), and (ii) uni-parametric linear programs (by replacing $Q(\theta)$ with the zero matrix in Problem \eqref{prob:QP}). To our knowledge, the only other works that propose solution strategies for Problem \eqref{prob:QP} under similar assumptions to those that are implied by our Assumption \ref{asm:sufficient} are those of \citet{ritter1962verfahren}, \citet{valiaho1985unified}, and \citet{jonker2001one}.

\section{Computational Results}\label{sec:implementation/results}

Using the Python programming language, we develop an implementation of the methods described in Algorithms \ref{alg:solveUpLCP}--\ref{alg:getExtremes} in which the processing of set $\mathcal{S}$ on line \ref{alg:solveUpLCP:line:while} of Algorithm \ref{alg:solveUpLCP} is performed in parallel, when desired. Interested readers may obtain the code from \url{https://github.com/Nadelgren/upLCP_solver}. The code is written so that the user may specify any one of three types of problems: (i) upLCP, (ii) upQP, or (iii) upLP. Problems given in the form of upQP or upLP are converted to upLCP prior to utilization of the methods proposed herein, but the computed solution is ultimately provided in the context of the originally presented problem. We note that solutions for problems in the form of upQP or upLP contain not only parametric values for the problem's original decision variables, but also for slack variables for each constraint and dual variables for each constraint and non-negativity restriction. 

As suggested in Section \ref{sec:algorithm}, our implementation utilizes the criss-cross method of \citet{den1993linear} to compute complementary bases on line \ref{alg:solveUpLCP:line:midpoint} of Algorithm \ref{alg:solveUpLCP}. Additionally, we employ the computer algebra system Pari/GP 2.14.0 \citep{PARI2} via the Python library CyPari2 for all symbolic algebra needed to compute and process (e.g., calculate roots, take derivatives, etc.) the polynomial functions that define each invariancy interval.

%\section{Computational Results}\label{sec:results}
We now present numerical results for two classes of problems. Instances of the first class are obtained from biobjective QPs having convex objectives that are scalarized using the weighted-sum approach and then reformulated as upLCP. We refer to these instances as \texttt{boQP} instances. Instances of the second class are obtained using some of the techniques outlined by \citet{illes2018generating}. We refer to these instances as \texttt{sufLCP} instances. All instances, together with a detailed description of the specific techniques use for their generation, are provided along with the code at the above referenced url. 

All tests were conducted on a machine running Linux Mint 20.0 and that had a 1.2 GHz Intel i3-1005G1 CPU with 12GB of RAM. When solving each instance of upLCP, we permitted the processing of set $\mathcal{S}$ to be executed in parallel, using four threads. In all, we generated five \texttt{boQP} instances for each value of $h \in \{50,75,100,125\}$ and five \texttt{sufLCP} instances for each value of $h \in \{50,75,100,125, 150, 175\}$. The CPU time required to solve each \texttt{boQP} instance, as well as the number of invariancy intervals present in the final partition for each of instance, is presented in Table \ref{tab:boqp}. Analogous data is provided for \texttt{sufLCP} instances in Table \ref{tab:suflcp}. 

\begin{table} \small
\centering
\caption{Results for \texttt{boQP} Instances}\label{tab:boqp}
\begin{tabular}{rrrr|rrrr}
Instance		& 		& CPU		& \# of		& Instance		& 		& CPU		& \# of\\
Size ($h$)	& \#		& Time (s)	& Intervals	& Size ($h$)		& \#		& Time (s)	& Intervals\\
\hline
\hline
50			& 1		& 18.61		& 15			& 100			& 1		& 915.22		& 33\\
			& 2 		& 24.90		& 13			&				& 2		& 1041.95	& 38\\
			& 3		& 51.05		& 30			& 				& 3		& 1377.16	& 39\\
			& 4		& 18.17		& 22			& 				& 4		& 806.51		& 27\\
			& 5		& 8.07		& 11			& 				& 5		& 1543.52	& 44\\
\hline
75			& 1		& 66.69		& 17			& 125			& 1		& 2813.07	& 38\\
			& 2 		& 185.42		& 27			& 				& 2		& 4463.14	& 51	\\
			& 3		& 327.48		& 23			& 				& 3		& 2137.46	& 45	\\
			& 4		& 266.87		& 26			& 				& 4		& 705.88		& 5	\\
			& 5		& 222.91		& 24			& 				& 5		& 1429.14	& 31	\\
\hline
\end{tabular}
\end{table}

\begin{table} \small
\centering
\caption{Results for \texttt{sufLCP} Instances}\label{tab:suflcp}
\begin{tabular}{rrrr|rrrr}
Instance		& 		& CPU		& \# of		& Instance		& 		& CPU		& \# of\\
Size ($h$)	& \#		& Time (s)	& Intervals	& Size ($h$)		& \#		& Time (s)	& Intervals\\
\hline
\hline
50			& 1		& 1.42		& 15			& 125			& 1		& 496.97		& 37\\
			& 2 		& 3.34		& 11			&				& 2		& 303.12		& 25\\
			& 3		& 5.73		& 17			& 				& 3		& 64.05		& 44\\
			& 4		& 2.81		& 10			& 				& 4		& 169.82		& 29\\
			& 5		& 3.94		& 6			& 				& 5		& 283.03		& 35\\
\hline
75			& 1		& 19.88		& 20			& 150			& 1		& 369.89		& 42\\
			& 2 		& 27.35		& 22			& 				& 2		& 814.69		& 48	\\
			& 3		& 39.17		& 52			& 				& 3		& 96.95		& 34\\
			& 4		& 6.23		& 9			& 				& 4		& 597.42		& 61	\\
			& 5		& 26.41		& 115		& 				& 5		& 442.07		& 37\\
\hline
100			& 1		& 80.94		& 23			& 175			& 1		& 1947.01	& 47\\
			& 2 		& 15.74		& 16			& 				& 2		& 2066.37	& 42\\
			& 3		& 250.15		& 37			& 				& 3		& 1083.17	& 42	\\
			& 4		& 22.27		& 14			& 				& 4		& 2228.88	& 54	\\
			& 5		& 41.50		& 18			& 				& 5		& 2062.75	& 40	\\
\hline
\end{tabular}
\end{table}

As we expect, we see that for both classes of problems the required CPU time and number of intervals present in the final solution increase as problem size increases.

\section{Conclusion}\label{sec:conclusion}

We have presented a new method for solving upLCP, upQP, and upLP that is capable of solving some classes of problems that were not able to be solved by any known uni-parametric methodology from the literature. Moreover, we have demonstrated empirically that our proposed technique can be used to solve relatively large instances in reasonable time. 
%Our implementation is made available for use by the research community.
%In the future, we aim to improve upon the methods presented herein and develop techniques that can be used to solve even larger instances.
%\begin{acknowledgements}
%The author would like to thank Jacob Adelgren for offering helpful advice and feedback regarding the implementation.
%\end{acknowledgements}



\bibliographystyle{plainnat}

\bibliography{/home/nate/Dropbox/citations}

%\begin{data}
%The code and datasets used in this work are available at \url{https://github.com/Nadelgren/upLCP_solver}.
%\end{data}
%
%\begin{declarations}
%The author did not receive support from any organization for the submitted work and has no relevant financial or non-financial interests to disclose.
%\end{declarations}

\end{document}